\documentclass[11pt]{article}
\usepackage{amssymb}
\newcommand{\be}{\begin{equation}}
\newcommand{\ee}{\end{equation}}
\newcommand{\bea}{\begin{eqnarray}}
\newcommand{\eea}{\end{eqnarray}}
\newcommand{\binom}[2]{{#1 \choose #2}}
\newcommand{\mbf}[1]{\mbox{\boldmath$#1$}}

\begin{document}

%\DeclareGraphicsExtensions{.jpg,.pdf,.png}

\title{Extension of the
Bernoulli and Eulerian Polynomials of Higher Order and
Vector Partition Function}
\author{Boris Y. Rubinstein
\\
Department of Mathematics, University of California, Davis,
\\One Shields Ave., Davis, CA 95616, U.S.A.}
\date{\today}

\maketitle
%%%%%%%%%%%%%%%%%%%%%%%%%%%%%%%%%%%%%%%%%%%%%%%%%%%%%%%%%%%%%%%%%%%%%
\begin{abstract}
Following the ideas of L. Carlitz we introduce a generalization
of the Bernoulli and Eulerian polynomials of higher order to
vectorial index and argument.
These polynomials are used for computation of the
{\it vector partition function} $W({\bf s},{\bf D})$, {\it i.e.},
a number of integer solutions to a linear system
${\bf x} \ge 0, \ \ {\bf D x} = {\bf s}$. It is shown that
$W({\bf s},{\bf D})$ can be expressed through
the {\it vector} Bernoulli polynomials of higher order.
\end{abstract}
%%%%%%%%%%%%%%%%%%%%%%%%%%%%%%%%%%%%%%%%%%%%%%%%%%%%%%%%%%%%%%%%%%%%%
%\newpage
\section{Introduction}
\label{intro}
The history of the Bernoulli polynomials $B_k(x)$
counts more than 250 years, as 
L. Euler first studied them for arbitrary values of the argument.
He introduced the Euler polynomials $E_k(x)$,
and relation between these functions was established in
the end of nineteenth century by L. Saalsch\"utz.  
The generalization of the Euler
polynomials to the so-called Eulerian polynomials
$H_k(x,\rho)$ was made by Frobenius  \cite{Frobenius}
several years later.

%In 1920s
N. N\"orlund \cite{NorlundMemo}
introduced the
Bernoulli $B^{(m)}_k(x|{\bf d}^m)$ 
and Euler polynomials $E^{(m)}_k(x|{\bf d}^m)$
of higher order adding 
$m$ parameters. 
Similar extension for the Eulerian polynomials 
$H^{(m)}_k(x,\mbf{\rho}^m|{\bf d}^m)$ was made by
L. Carlitz in \cite{Carlitz1960a}.
%The polynomials of higher order, their 
%recursive and symmetry properties and some useful
%formulas are presented in Section \ref{PolyScal1}.

The Bernoulli and Eulerian polynomials of higher order
appear to be useful for the description of the
{\it restricted partition function} $W(s,{\bf d}^m)$,
which is a number of integer nonnegative solutions of 
Diophantine equation ${\bf d}^m\cdot{\bf x}=s$
(see \cite{Rama03}).
The author showed in \cite{Rub04}
that $W(s,{\bf d}^m)$ may be written as a finite sum of the 
Bernoulli polynomial of higher order multiplied by
{\it prime circulator} introduced by A. Cayley
(see \cite{Dickson1923}).
The short derivation of this result
is given in Section \ref{ScalarPartitionFunction}.

Carlitz suggested in \cite{Carlitz1960b}
another extension of the Bernoulli
$B_{\bf k}({\bf x})$ and 
Eulerian polynomials $H_{\bf k}({\bf x},\rho)$ 
to vectorial arguments and
indices. % which are described in Section \ref{PolyVect0}.

In this work we combine both abovementioned approaches
to introduce a new class of polynomials, 
which we call {\it vector} Bernoulli
$B^{(l,m)}_{\bf k}({\bf x}|{\bf D})$
and Eulerian 
$H^{(l,m)}_{\bf k}({\bf x},\mbf{\rho}|{\bf D})$
polynomials of higher order. 
We define the generating functions for these
polynomials, find the recursive and symmetry relations,
and find a new relation between
these polynomials.
Using the vector polynomials of higher order in
Section \ref{VectorPartitionFunction} we find an explicit
formula for the {\it vector partition function} $W({\bf s},{\bf D})$,
defined as a number of integer nonnegative solutions to the linear system
${\bf x} \ge 0, \ \ {\bf D x} = {\bf s}$,
where ${\bf D}$ denotes a nonnegative non-degenerate
integer %$(l \times m)$ 
matrix.
It appears that similar to the scalar case 
the vector partition function
may be written as a finite sum of the
vector Bernoulli polynomials of higher order multiplied by
prime circulators of vector index and argument. 
The solution gives the vector partition function in every
chamber of the system, as well as it determines the shape of
each individual chamber.

%%%%%%%%%%%%%%%%%%%%%%%%%%%%%%%%%%%%%%%%%%%%%%%%%%%%%%%%%%%%%%%%%%%%%%%%%%%%%%%%%

%%%%%%%%%%%%%%%%%%%%%%%%%%%%%%%%%%%%%%%%%%%%%%%%%%%%%%%%%%%%%%%%%%%%%%%%%%%%%%%%%
\section{Bernoulli and Eulerian Polynomials and Their Generalizations}
\label{BEpoly}
In the Section \ref{PolyScal0}
the definition and main properties of
the (regular) Bernoulli $B_k(x)$
and Eulerian $H_k(x,\rho)$ polynomials are presented. 
Their generalization 
to the polynomials of
higher orders $B^{(m)}_k(x|{\bf d}^m)$ and
$H^{(m)}_k(x,\mbf{\rho}^m|{\bf d}^m)$,
recursive and symmetry properties and some useful
formulas are presented in Section \ref{PolyScal1}. 
Extension of the regular polynomials
to vectorial arguments and
indices $B_{\bf k}({\bf x})$ and $H_{\bf k}({\bf x},\rho)$
described in the Section \ref{PolyVect0}.
Finally, in the Section \ref{PolyVect1} combining both extensions
we introduce a new class of the polynomials
$B^{(l,m)}_{\bf k}({\bf x}|{\bf D})$ and
$H^{(l,m)}_{\bf k}({\bf x},\mbf{\rho}|{\bf D})$ and
consider their properties.

%%%%%%%%%%%%%%%%%%%%%%%%%%%%%%%%%%%%%%%%%%%%%%%%%%%%%%%%%%%%%%%%%%%%%%%%%%%%%%%%%
\subsection{Bernoulli and Eulerian Polynomials}
\label{PolyScal0}
Define the Bernoulli polynomial $B_k(x)$
through the generating function
$$
%\be
e^{xt}\frac{t}{e^t-1} = \sum_{k=0}^{\infty} B_k(x) \frac{t^k}{k!}.
$$
%\label{Bern0GF}
%\ee
The Bernoulli numbers are defined as $B_k = B_k(0)$.
It is known that the {\it multiplication formula} holds for these polynomials:
$$
%\be
\sum_{r=0}^{p-1} B_k\left(x+\frac{r}{p} \right)=p^{1-k} B_k(px);
$$
%\label{Bern0mult}
%\ee
this relation can be used as an alternative definition
\cite {Nielsen1923} of the Bernoulli polynomials.
Another remarkable property of the Bernoulli polynomials is
that they can be presented as a symbolic power
$$
%\be
B_k(x) = (B + x)^k = \sum_{i=0}^k {\binom k i} B_i x^{k-i},
$$
%\label{Bern0symb}
%\ee
where after the expansion the powers of $B$ are transformed into the indices.

The Eulerian polynomials were considered by Frobenius in \cite{Frobenius}, their
generating function reads
$$
%\be
e^{xt}\frac{1-\rho}{e^t-\rho} =
\sum_{k=0}^{\infty} H_k(x,\rho) \frac{t^k}{k!}, \ \
(\rho \ne 1).
$$
%\label{Euler0GF}
%\ee
For $\rho=-1$ the Eulerian polynomials reduce to the Euler polynomials
$E_k(x)$. The former also admits the symbolic form
$$
%\be
H_k(x,\rho) = (H(\rho) + x)^k,
$$
%\label{Euler0symb}
%\ee
where $H_k(\rho) \equiv H_k(0,\rho)$ are the Eulerian numbers.
For $\rho^p = 1$ we have an additional relation \cite{Carlitz1960a}
between the Bernoulli and Eulerian polynomials
\be
p^{k-1} \sum_{r=0}^{p-1} \rho^{-r} B_k\left(x+\frac{r}{p} \right)=
\frac{k\rho}{1-\rho} H_{k-1}(px,\rho).
\label{BernEuler0Rel}
\ee

%%%%%%%%%%%%%%%%%%%%%%%%%%%%%%%%%%%%%%%%%%%%%%%%%%%%%%%%%%%%%%%%%%%%%%%%%%%%%%%%%
\subsection{Bernoulli and Eulerian Polynomials of Higher Order}
\label{PolyScal1}
N\"orlund \cite{NorlundMemo} introduced the Bernoulli
polynomials of higher order $B^{(m)}(x|{\bf d}^m)$
adding $m$ parameters ${\bf d}^m=\{d_1,d_2,\ldots,d_m\}$ by means of the
generating function:
$$
%\be
e^{xt}\prod_{i=1}^m \frac{d_i t}{e^{d_i t}-1} =
\sum_{k=0}^{\infty} B^{(m)}_k(x|{\bf d}^m) \frac{t^k}{k!}.
$$
%\label{Bern1GF}
%\ee
The regular Bernoulli polynomials are
expressed through the relation
$$
B^{(1)}_k(x|d) = d^k B_k(x).
$$
N\"orlund also found the multiplication formulas
$$
%\be
\sum_{r=0}^{p-1}
B^{(m)}_k\left(x+\frac{r}{p}\sigma({\bf d}^m)|{\bf d}^m\right)
=p^{m-k} B^{(m)}_k(px|{\bf d}^m), \ \ \
\sigma({\bf d}^{m}) = \sum_{i=1}^m d_i,
$$
%\label{Bern1mult1}
%\ee
$$
%\be
\sum_{r_i=0}^{p_i-1}
B^{(m)}_k(x+{\bf r}\cdot{\bf d}^m|\{p_1d_1,\ldots,p_md_m\})=
\pi({\bf p}^m) B^{(m)}_k(x|{\bf d}^m), \ \ \
\pi({\bf p}^m) = \prod_{i=1}^m p_i,
$$
%\label{Bern1mult2}
%\ee
where ${\bf r}\cdot{\bf d}^m = \sum_{i=1}^m r_i d_i$ denotes the scalar product.
N\"orlund gave the symbolic expression 
% for the polynomials
$$
%\be
B^{(m)}_k(x|{\bf d}^m) =
(x+ \sum_{i=1}^m d_i \;{}^i\! B)^k,
$$
%\label{Bern1symb}
%\ee
where again after the expansion powers of ${}^i\! B$ are converted into the
indices of the Bernoulli numbers.
This relation implies a binomial formula:
\be
\sum_{k=0}^{n} \binom{n}{k}
B^{(m)}_k(x|{\bf d}^m) B^{(l)}_{n-k}(y|{\bf d}^l) =
B^{(m+l)}_n(x+y|{\bf d}^{m+l}).
\label{Bern1binom}
\ee
The following symmetry relation is valid
\be
B^{(m)}_k(x|{\bf d}^m) =
(-1)^k B^{(m)}_{k}(-x-\sigma({\bf d}^{m})|{\bf d}^{m}),
\label{Bern1symm}
\ee
and the recursive relation holds for the Bernoulli polynomials
\be
B^{(m)}_k(x+d_m|{\bf d}^m)-B^{(m)}_k(x|{\bf d}^m) =
k d_m B^{(m-1)}_{k-1}(x|{\bf d}^{m-1}).
\label{Bern1recurs}
\ee
%N\"orlund also found similar generalization of
%the Euler polynomials.

The Eulerian polynomials of higher order
$H^{(m)}(x, \mbf{\rho}^m|{\bf d}^m)$ were considered by Carlitz in
\cite{Carlitz1960a}. Their generating function reads
$$
%\be
e^{xt}\prod_{i=1}^m\frac{1-\rho_i}{e^{d_i t}-\rho_i} =
\sum_{k=0}^{\infty}
H^{(m)}_k(x, \mbf{\rho}^m|{\bf d}^m)\frac{t^k}{k!},
$$
%\label{Euler1GF}
%\ee
%Carlitz found the symbolic form for these functions
and they admit the following symbolic notation:
$$
%\be
H^{(m)}_k(x, \mbf{\rho}^m|{\bf d}^m)=
(x+ \sum_{i=1}^m d_i H(\rho_i))^k.
$$
%\label{Euler1symb}
%\ee
When $\rho_j^{p_j} = 1, \ \rho_j \ne 1, \ j=1,2,\ldots,m$ we
have the analog of
(\ref{BernEuler0Rel})
\bea
&&
\frac{1}{\pi({\bf p}^m)}
\sum_{r_j=0}^{p_j-1}
\rho_j^{-r_j} B^{(m)}_k(x+{\bf r}\cdot{\bf d}^m|
\{p_1d_1,p_2d_2,\ldots,p_md_m\})
\nonumber \\
&=&
\frac{k!}{(k-m)!} \left(\prod_{i=1}^m \frac{\rho_i d_i}{1-\rho_i}\right)
H^{(m)}_{k-m}(x,\mbf{\rho}^m|{\bf d}^m).
\label{BernEuler1Rel}
\eea

%%%%%%%%%%%%%%%%%%%%%%%%%%%%%%%%%%%%%%%%%%%%%%%%%%%%%%%%%%%%%%%%%%%%%%%%%%%%%%%%%
\subsection{Vector Bernoulli and Eulerian Polynomials}
\label{PolyVect0}
Carlitz %in \cite{Carlitz1960b}
suggested another extension of the
Bernoulli and Eulerian polynomials.
%, namely, to the vector argument.
He introduced \cite{Carlitz1960b} the Bernoulli polynomials $B_{\bf k}({\bf x})$
of the vector argument
${\bf x} = \{x_1,\ldots, x_l\}$ and vector index
${\bf k} = \{k_1,\ldots, k_l\}$ through the generating function
$$
%\be
e^{\bf x \cdot t} \frac{\sum_{i=1}^l t_i}{e^{\sum_{i=1}^l t_i}-1}=
\sum_{\bf k}B_{\bf k}({\bf x})
%\prod_{i=1}^l\frac{t_i^{k^i}}{k_i!},
\frac{{\bf t^k}}{{\bf k}!},
$$
%\label{Bern0vGF}
%\ee
where the summation is over all $k_1,\ldots,k_l$ and
${\bf x \cdot t} = \sum_{i=1}^l x_i t_i$ is a scalar product.
Here and below we use the notation
$$
{\bf t^k} = \prod_{i=1}^l t_i^{k_i}, \ \
{\bf k}! = \prod_{i=1}^l k_i!, \ \
\binom{\bf n}{\bf k}=\prod_{i=1}^l \binom{n_i}{k_i}.
$$
The multiplication theorem is established in the form
$$
%\be
\sum_{r=0}^{p-1} B_{\bf k}\left({\bf x}+\frac{r}{p} \right) =
p^{1-|{\bf k}|} B_{\bf k}(p{\bf x}),\ \
|{\bf k}| = \sum_{i=1}^l k_i.
$$
%\label{Bern0vmult}
%\ee

The Eulerian polynomials $H_{\bf k}({\bf x},\rho)$
of the vector argument and index are defined as
$$
%\be
e^{\bf x \cdot t} \frac{1-\rho}{e^{\sum_{i=1}^l t_i}-\rho}=
\sum_{{\bf k}=0}^{\infty}H_{\bf k}({\bf x},\rho)
%\prod_{i=1}^l\frac{t_i^{k^i}}{k_i!},
\frac{{\bf t^k}}{{\bf k}!},
$$
%\label{Euler0vGF}
%\ee
and the analog of (\ref{BernEuler0Rel}) for $\rho^p=1$ reads:
$$
%\be
p^{|{\bf k}|-1}\sum_{r=0}^{p-1}
\rho^{-r} B_{\bf k}\left({\bf x}+\frac{r}{p}\right)=
\frac{\rho}{1-\rho} \sum_{i=1}^l k_i H_{{\bf k}-{\bf 1}^{(i)}}(p{\bf x},\rho),
$$
%\label{BernEulerVect0Rel}
%\ee
where
${\bf 1}^{(i)}$ is vector with the components
${\bf 1}^{(i)}_m = \delta_{im}$, and $\delta_{im}$ is the Kroneker delta.

Here and below we will use the notion of the {\it vector
Bernoulli $B_{\bf k}({\bf x})$ and Eulerian $H_{\bf k}({\bf x},\rho)$
polynomials} remembering
that they are scalar functions of
their arguments.

%%%%%%%%%%%%%%%%%%%%%%%%%%%%%%%%%%%%%%%%%%%%%%%%%%%%%%%%%%%%%%%%%%%%%%%%%%%%%%%%%
\subsection{Vector Bernoulli and Eulerian Polynomials of higher order}
\label{PolyVect1}
In this work we make one more step in generalization of the
Bernoulli and Eulerian
polynomials which naturally arises from the approaches
discussed above. We consider {\it vector
Bernoulli $B^{(l,m)}_{\bf k}({\bf x}|{\bf D}^m)$ and
Eulerian $H^{(l,m)}_{\bf k}({\bf x},\mbf{\rho}^m|{\bf D}^m)$
polynomials of higher order}, where ${\bf D}^m$ denotes $(l \times m)$ matrix.
%with integer nonnegative entries.

The vector Bernoulli polynomials of higher order are defined by means of
the generating function:
$$
%\be
e^{\bf x \cdot t} \prod_{j=1}^{m}
\frac{\sum_{i=1}^l t_i D_{ij}}{e^{\sum_{i=1}^l t_i D_{ij}}-1}=
\sum_{{\bf k}}B^{(l,m)}_{\bf k}({\bf x}|{\bf D}^m)
%\prod_{i=1}^l\frac{t_i^{k^i}}{k_i!}\;.
\frac{{\bf t^k}}{{\bf k}!}\;.
$$
%\label{Bern1vGF}
%\ee
Denoting the columns of the matrix ${\bf D}^m$ as
$\{{\bf c}_1, \ldots, {\bf c}_m$\}, we can rewrite the above definition as
\be
e^{\bf x \cdot t} \prod_{j=1}^{m}
\frac{{\bf t \cdot c}_j}{e^{{\bf t \cdot c}_j}-1}=
\sum_{{\bf k}}B^{(l,m)}_{\bf k}({\bf x}|{\bf D}^m)
%\prod_{i=1}^l\frac{t_i^{k^i}}{k_i!}\;.
\frac{{\bf t^k}}{{\bf k}!}\;.
\label{Bern1vGF}
\ee
The analog of the recursive relation (\ref{Bern1recurs}) reads
$$
%\be
B^{(l,m)}_{\bf k}({\bf x}+{\bf c}_m|{\bf D}^m) -
B^{(l,m)}_{\bf k}({\bf x}|{\bf D}^m) =
\sum_{i=1}^l k_i D_{im}
B^{(l,m-1)}_{{\bf k}-{\bf 1}^{(i)}}({\bf x}|{\bf D}^{m-1})\;.
$$
%\label{Bern1vrecurs}
%\ee
%where
%${\bf 1}^{(i)}$ is vector with the components
%${\bf 1}^{(i)}_m = \delta_{im}$, and $\delta_{im}$ is the Kroneker delta.
One can check the validity of an analog of symmetry relation
(\ref{Bern1symm})
$$
%\be
B^{(l,m)}_{\bf k}({\bf x}|{\bf D}^m) =
(-1)^{|{\bf k}|}
B^{(l,m)}_{\bf k}(-{\bf x}-\mbf{\sigma}({\bf D}^m)|{\bf D}^m), \ \ \
\mbf{\sigma}({\bf D}^m) = \sum_{j=1}^m {\bf c}_j.
$$
%\label{Bern1vsymm}
%\ee
Using (\ref{Bern1vGF}) we obtain the vector generalization of the 
binomial formula (\ref{Bern1binom})
\be
\sum_{\bf k}^{\bf n} \binom{\bf n}{\bf k}
B^{(l,m_1)}_{\bf k}({\bf x}|{\bf D}^{m_1})
B^{(l,m_2)}_{{\bf n-k}}({\bf y}|{\bf D}^{m_2}) =
B^{(l,m_1+m_2)}_{\bf n}({\bf x}+{\bf y}|{\bf D}^{m_1+m_2}).
\label{Bern1vbinom}
\ee

The vector Eulerian polynomials of higher order have the generating function of
the form:
$$
%\be
e^{\bf x \cdot t} \prod_{j=1}^{m}
\frac{1-\rho_j}{e^{{\bf t \cdot c}_j}-\rho_j}=
\sum_{{\bf k}}H^{(l,m)}_{\bf k}({\bf x},\mbf{\rho}^{m}|{\bf D}^m)
%\prod_{i=1}^l\frac{t_i^{k^i}}{k_i!}\;,
\frac{{\bf t^k}}{{\bf k}!}\;,
$$
%\label{Euler1vGF}
%\ee
where $\mbf{\rho}^{m} =\{\rho_1,\ldots,\rho_m\}$.
Consider a relation between the vector Bernoulli and Eulerian polynomials
of higher order similar to (\ref{BernEuler1Rel}). For $p_j$ such that
$\rho_j^{p_j}=1$ we have the following expression:
\bea
&&
\sum_{\bf k}\sum_{r_j=0}^{p_j-1}\rho_j^{-r_j}
B^{(l,m)}_{\bf k}({\bf x}+\sum_{j=1}^m r_j {\bf c}_j|
\{p_1{\bf c}_1,\ldots,p_m{\bf c}_m\})
\frac{{\bf t^k}}{{\bf k}!}
\nonumber \\
&=&
\prod_{j=1}^m \sum_{r_j=0}^{p_j-1}\rho_j^{-r_j}
\exp[({\bf x}+\sum_{j=1}^m r_j {\bf c}_j)\cdot {\bf t}]
\left(\frac{p_j {\bf c}_j\cdot {\bf t}}
{e^{p_j {\bf c}_j\cdot {\bf t}}-1} \right)
\nonumber \\
&=&
\left(\prod_{j=1}^m
\frac{p_j\rho_j}{1-\rho_j} ({\bf c}_j\cdot {\bf t})
\right)
e^{{\bf x}\cdot {\bf t}}
\prod_{j=1}^m \frac{1-\rho_j}{e^{{\bf c}_j\cdot {\bf t}}-\rho_j}
\nonumber \\
&=&
\pi({\bf p}^m)
\left(
\prod_{j=1}^m \frac{\rho_j}{1-\rho_j}({\bf c}_j\cdot {\bf t})
\right)
\sum_{\bf k}
H^{(l,m)}_{\bf k}({\bf x},\mbf{\rho}^{m}|{\bf D}^m)
\frac{{\bf t^k}}{{\bf k}!}.
\label{BernEuler1vRel}
\eea

%%%%%%%%%%%%%%%%%%%%%%%%%%%%%%%%%%%%%%%%%%%%%%%%%%%%%%%%%%%%%%%%%%%%%%%%%%%%%%%%%

%%%%%%%%%%%%%%%%%%%%%%%%%%%%%%%%%%%%%%%%%%%%%%%%%%%%%%%%%%%%%%%%%%%%%%%%%%%%%%%%%
\section{Restricted Partition Function}
\label{ScalarPartitionFunction}
The restricted
partition function $W(s,{\bf d}^m) \equiv W(s,\{d_1,d_2,\ldots,d_m\})$ is a
number of partitions of an integer $s$ into positive integers
$\{d_1,d_2,\ldots,d_m\}$,
each not greater than $s$. The generating function for $W(s,{\bf d}^m)$
has a form
$$
%\be
\prod_{i=1}^m\frac{1}{1-t^{d_{i}}}
 =\sum_{s=0}^{\infty} W(s,{\bf d}^m)\;t^s\;,
$$
%\label{WGF}
%\ee
where $W(s,{\bf d}^m)$ satisfies the basic recursive relation
\be
W(s,{\bf d}^m) - W(s-d_m,{\bf d}^m) = W(s,{\bf d}^{m-1})\;.
\label{Wrecurs}
\ee
Note a similarity of (\ref{Wrecurs}) to the recursive relation
(\ref{Bern1recurs}) for the Bernoulli polynomials of higher order.
Sylvester found a symmetry property of the partition
function:
\be
W(s,{\bf d}^m) = (-1)^{m-1} W(-s-\sigma({\bf d}^m),{\bf d}^m), \ \ \
\sigma({\bf d}^m)= \sum_{i=1}^m d_i.
\label{Wsymm}
\ee
He proved \cite{Sylv2} a statement about splitting of the partition
function into periodic and non-periodic parts and showed that the
restricted partition function may be presented as a sum of "waves", which
we call the {\em Sylvester waves}
\be
W(s,{\bf d}^m) = \sum_{j=1} W_j(s,{\bf d}^m)\;,
\label{SylvWavesExpand}
\ee
where summation runs over all distinct factors
of the elements of the set ${\bf d}^m$.
The wave $W_j(s,{\bf d}^m)$ is a quasipolynomial in $s$
closely related to prime roots $\rho_j$ of unity.
The wave
$W_j(s,{\bf d}^m)$ is a coefficient of
${t}^{-1}$ in the series expansion in ascending powers of $t$ of
the generator
\be
F_j^m(s,t)=
\sum_{\rho_j} \frac{\rho_j^{-s} e^{st}}{\prod_{k=1}^{m}
        \left(1-\rho_j^{d_k} e^{-d_k t}\right)}\;.
\label{generatorWj}
\ee
The summation is made over all prime roots of unity
$\rho_j=\exp(2\pi i n/j)$ for $n$ relatively prime to $j$
(including unity) and smaller than $j$. 
It is easy to check by straightforward calculation that 
the recursive relation
\be
F_j^m(s,t)-F_j^m(s-d_m,t)=F_j^{m-1}(s,t)
\label{Frecurs}
\ee
holds for any generator $F_j^m(s,t)$,
impliyng the validity of (\ref{Wrecurs})
for each Sylverster wave $W_j(s,{\bf d}^m)$.
The generator satisfies the following symmetry property
$$
F_j^m(s,t)=(-1)^m F_j^m(-s-\sigma({\bf d}^m),-t),
$$
which implies the validity of (\ref{Wsymm})
for the residue of $F_j^m(s,t)$.

In \cite{Rama03} the
explicit expression for the Sylvester wave of the arbitrary period is given
through the Bernoulli and Eulerian polynomials of
higher order.
% It was also proved that both the recursion (\ref{Wrecurs})
% and the symmetry relation (\ref{Wsymm})
% appear to be valid for
% every Sylvester wave entering (\ref{SylvWavesExpand}).
Using (\ref{BernEuler1Rel})
it was shown in \cite{Rub04} that it is possible to
express the Sylvester wave as a finite sum of the Bernoulli polynomials of
higher order only.
Here we present a short derivation of this result. 

Assuming that the vector ${\bf d}^m$ has
$\omega_j$ components divisible
by $j$,
sort the elements of ${\bf d}^m$ in such way that
the $j$-divisible integers come first. The generator
(\ref{generatorWj}) can be written as a product
$$
%\be
F_j(s,t) = \sum_{\rho_j}
\frac{e^{st}}{\prod_{i=1}^{\omega_j} (1-e^{-d_it})} \times
\frac{\rho_j^{-s}}{\prod_{i=\omega_j+1}^m (1-\rho_j^{d_i}e^{-d_i t})}\;.
$$
%\label{generatorWj_product}
%\ee
%Note that the
%residue of the generator (\ref{generatorWj}) is equal to the coefficient
%of $t^{m-1}$ of the expression 
Consider a modified generator
$\tilde{F}_j^m(s,t)=\pi({\bf d}^m)t^m F_j^m(s,t),$
% divided by $\pi({\bf d}^m)$.
% for which we have a chain of transformations below. 
for which using the notation 
$
\rho^{{\bf d}^m} = \{\rho^{d_1},\rho^{d_2},\ldots,\rho^{d_m}\}
$
we have
\bea
\tilde{F}_j^m(s,t) &=&
t^{m-\omega_j} \sum_{\rho_j}
\rho_j^{-s} \prod_{i=1}^{\omega_j} \frac{d_it}{e^{d_i t}-1}\cdot
e^{(s+\sigma({\bf d}^m))t}\prod_{i=\omega_j+1}^m
\frac{d_i}{e^{d_i t}-\rho_j^{d_i}}
\nonumber \\
&=&
t^{m-\omega_j} \sum_{\rho_j}
\rho_j^{-s-\sigma({\bf d}^{m-\omega_j})}
\prod_{i=1}^{\omega_j} \frac{d_it}{e^{d_i t}-1}\cdot
e^{(s+\sigma({\bf d}^m))t}\prod_{i=\omega_j+1}^m
\frac{1-\rho_j^{d_i}}{e^{d_i t}-\rho_j^{d_i}}\cdot
\frac{d_i \rho_j^{d_i}}{1-\rho_j^{d_i}}
\nonumber \\
&=&
t^{m-\omega_j} \sum_{\rho_j}
\rho_j^{-s-\sigma({\bf d}^{m-\omega_j})}
\cdot
\sum_{n_1=0}^{\infty} 
B_{n_1}^{(\omega_j)}(0|{\bf d}^{\omega_j})\frac{t^{n_1}}{n_1!}
\label{Wj1} \\
&\times&
\prod_{i=\omega_j+1}^m \frac{d_i \rho_j^{d_i}}{1-\rho_j^{d_i}}\cdot
\sum_{n_2=0}^{\infty}
H_{n_2}^{(m-\omega_j)}(s+\sigma({\bf d}^m),
\rho_j^{{\bf d}^{m-\omega_j}}|{\bf d}^{m-\omega_j})
\frac{t^{n_2}}{n_2!}.
\nonumber
\eea
Now we employ (\ref{BernEuler1Rel}), noting that
for $\omega+1\le i \le m$ all $p_i=j$, to obtain
\bea
&&\prod_{i=\omega_j+1}^m \frac{d_i \rho_j^{d_i}}{1-\rho_j^{d_i}}\cdot
\sum_{n_2=0}^{\infty}
H_{n_2}^{(m-\omega_j)}(s+\sigma({\bf d}^m),
\rho_j^{{\bf d}^{m-\omega_j}}|{\bf d}^{m-\omega_j})
\frac{t^{n_2}}{n_2!}
\nonumber \\
&=&
\sum_{n_2=0}^{\infty}
\frac{t^{n_2}}{(m-\omega_j+n_2)!j^{m-\omega_j}}
\sum_{r_i=0}^{j-1}
\rho_j^{-r_id_i} B_{m-\omega_j+n_2}^{(m-\omega_j)}
(s+\sigma({\bf d}^m)+{\bf r}\cdot {\bf d}^{m-\omega_j}
|j{\bf d}^{m-\omega_j}).
\nonumber
\eea
Inserting the above expression into (\ref{Wj1}) and using
(\ref{Bern1binom}) we obtain
\be
\tilde{F}_j^m(s,t) =
j^{-(m-\omega_j)}
\sum_{\rho_j}\sum_{r_i=0}^{j-1}
\rho_j^{-s-{\bf r}\cdot {\bf d}^{m-\omega_j}-\sigma({\bf d}^m)}
%\sum_{i=\omega_j+1}^m (r_i+1)d_i}
\sum_{n} B_n^{(m)} (s+\sigma({\bf d}^m)+
{\bf r}\cdot {\bf d}^{m-\omega_j}|{\bf d}^m_j)
\frac{t^{n}}{n!},
\label{Wj2}
\ee
where we use a shorthand notation ${\bf d}^m_j$ for
{\it $j$-modified} set of summands
defined as union of subset ${\bf d}^{\omega_j}$ of
summands divisible by $j$ and
the remaining part ${\bf d}^{m-\omega_j}$ multiplied by $j$
$$
{\bf d}^m_j = {\bf d}^{\omega_j} \cup j{\bf d}^{m-\omega_j} =
\{d_1,\ldots,d_{\omega_j},jd_{\omega_j+1},\ldots,jd_m\},
$$
so that ${\bf d}^m_j$ is divisible by $j$.
From (\ref{Wj2}) the expression for ${F}_j^m(s,t)$ follows:
\be
{F}_j^m(s,t) =
\frac{1}{\pi({\bf d}^m)j^{m-\omega_j}}
\sum_{\rho_j}\sum_{r_i=0}^{j-1}
\rho_j^{-s-{\bf r}\cdot {\bf d}^{m-\omega_j}-\sigma({\bf d}^m)}
%\sum_{i=\omega_j+1}^m (r_i+1)d_i}
\sum_{n} B_n^{(m)} (s+\sigma({\bf d}^m)+
{\bf r}\cdot {\bf d}^{m-\omega_j}|{\bf d}^m_j)
\frac{t^{n-m}}{n!}.
\label{Wj3}
\ee
Setting in (\ref{Wj3}) $n=m-1$ 
% and dividing the corresponding coefficient by $\pi({\bf d}^m)$
we arrive at the expression for the Sylvester wave:
$$
%\be
W_j(s,{\bf d}^m)  =
\frac{1}{(m-1)! \; \pi({\bf d}^m)\;j^{m-\omega_j}}
\sum_{r_i=0}^{j-1}
B_{m-1}^{(m)} (s+\sigma({\bf d}^m)+
{\bf r}\cdot {\bf d}^{m-\omega_j}|{\bf d}^m_j)
\sum_{\rho_j}
\rho_j^{-s-{\bf r}\cdot
{\bf d}^{m-\omega_j}-\sigma({\bf d}^m)}.
%\sum_{i=\omega_j+1}^m (r_i+1)d_i}.
$$
%\label{Wjfinal}
%\ee
Introduce a notation
$$%\be
\Psi_j(s) = % \frac{1}{j}
\sum_{\rho_j} \rho_j^s
%\label{Psibar_def}
$$
for the {\em prime radical circulator} (see \cite{Dickson1923}).
For prime $j$ it is given by
$$
%\be
\Psi_j(s) =
	\left\{ \begin{array}{ll}
         \phi(j), & \mbox{$s=0 \pmod{j}$}, \\
	 \mu(j), & \mbox{$s\ne 0 \pmod{j}$},
         \end{array}\right.
$$
%\label{primecirc}
%\ee
where $\phi(j)$ and $\mu(j)$ denote Euler totient and M\"obius functions.
Considering $j$ as a product of powers of
distinct prime factors
$$
j = \prod_{k} p_{k}^{\alpha_{k}},
$$
one may easily check that for integer values of $s$
$$
%\be
\Psi_j(s) =
\prod_{k} p_{k}^{\alpha_{k}-1}
\Psi_{p_{k}}\left(\frac{s}{p_{k}^{\alpha_{k}-1}}\right),
$$
%\label{gencirc}
%\ee
where $\Psi_k(s) = 0$ for non-integer values of $s$.

Noting that for $j$-modified set ${\bf d}^m_j$ we have
$$
\pi({\bf d}^m_j)=j^{m-\omega_j}\pi({\bf d}^m), \ \ \
\sigma({\bf d}^m_j)=\sigma({\bf d}^m)+(j-1)\sum_{i=\omega_j+1}^m d_i,
$$
and using the prime circulator notation 
we can write the Sylvester wave in a form
$$
%\be
W_j(s,{\bf d}^m)  =
\frac{1}{(m-1)! \; \pi({\bf d}^m_j)}
\sum_{r_i=0}^{j-1}
B_{m-1}^{(m)} (s+\sigma({\bf d}^m_j)-
{\bf r}\cdot {\bf d}^{m-\omega_j}
|{\bf d}^m_j)
\Psi_j(s-{\bf r}\cdot {\bf d}^{m-\omega_j}).
$$
%\label{Wjfinal1}
%\ee
The polynomial part of the partition function corresponds to 
$j=1$ and equals to
\be
W_1(s,{\bf d}^m) =
\frac{1}{(m-1)! \; \pi({\bf d}^m)}
B_{m-1}^{(m)} (s+\sigma({\bf d}^m)|{\bf d}^m).
\label{W1}
\ee
The polynomial part for the $j$-modified set ${\bf d}^m_j$ for $j>1$ reads
\be
W_1(s,{\bf d}^m_j) =
\frac{1}{(m-1)! \; \pi({\bf d}^m_j)}
B_{m-1}^{(m)} (s+\sigma({\bf d}^m_j)|{\bf d}^m_j).
\label{W1j}
\ee
Thus, the Sylvester wave for $j>1$ can be written as a linear superposition
of the polynomial part of the $j$-modified set multiplied by the 
corresponding prime circulator:
\be
W_j(s,{\bf d}^m) = \sum_{r_i=0}^{j-1}
W_1(s-{\bf r}\cdot {\bf d}^{m-\omega_j}
%\!\!\!\sum_{i=\omega_j+1}^m \!\!\!r_id_i
,{\bf d}^m_j)
\Psi_j(s-{\bf r}\cdot {\bf d}^{m-\omega_j}).
%\!\!\!\sum_{i=\omega_j+1}^m \!\!\!r_id_i).
\label{WjthroughW1j}
\ee
It is easy to see that for $j=1$ one has
${\bf d}^m_1 \equiv {\bf d}^m$, $\Psi_1(s)=1$ and
$\omega_1=m$, so that summation signs disappear, and
(\ref{WjthroughW1j}) reduces to (\ref{W1}).

%%%%%%%%%%%%%%%%%%%%%%%%%%%%%%%%%%%%%%%%%%%%%%%%%%%%%%%%%%%%%%%%%%%%%%%%%%%%%%%%%

%%%%%%%%%%%%%%%%%%%%%%%%%%%%%%%%%%%%%%%%%%%%%%%%%%%%%%%%%%%%%%%%%%%%%%%%%%%%%%%%%
\section{Restricted Vector Partition Function}
\label{VectorPartitionFunction}
Consider a function $W({\bf s},{\bf D}^m)$ counting the number of integer
nonnegative
solutions ${\bf x} \ge 0$
to the linear system ${\bf D}^m \cdot {\bf x} = {\bf s}$, where
${\bf D}^m$ is a nonnegative integer $(l \times m)$ matrix.
The function $W({\bf s},{\bf D}^m)$ is called {\it vector partition function} as
it is natural generalization of the restricted partition function to the 
vector argument.

The generating function for the vector partition function reads
\be
\prod_{i=1}^m \frac{1}{1-{\bf t}^{{\bf c}_i}} =
\sum_{{\bf s}} W({\bf s},{\bf D}^m) {\bf t^s} =
\sum_{{\bf s}} W({\bf s},\{{\bf c}_1,\ldots,{\bf c}_m\}) {\bf t^s}\;,
\label{WvectGF}
\ee
where ${\bf c}_i$ denotes the $i$-th column of the matrix ${\bf D}^m$.
It is easy to see that the vector analog of the
recursive relation (\ref{Wrecurs}) holds
\be
W({\bf s},{\bf D}^m) -
W({\bf s} - {\bf c}_m,{\bf D}^m) =
W({\bf s},{\bf D}^{m-1}).
\label{Wvectrecurs}
\ee
The symmetry property established in \cite{Beck2002} is just a vector
generalization of (\ref{Wsymm}):
\be
W({\bf s},{\bf D}^m) = (-1)^{m - \mbox{\small rank}\; {\bf D}^m}
W(-{\bf s}-\mbf{\sigma}({\bf D}^m),{\bf D}^m), \ \
\mbf{\sigma}({\bf D}^m) = \sum_{i=1}^m {\bf c}_i.
\label{Wvectsymm}
\ee
%For $m=1, l \ge 1$ we have ${\bf D}^1={\bf c}$, and
%\be
%W({\bf s},{\bf c}) = \Psi_{\bf c}({\bf s})=\prod_{k=1}^l \Psi_{c_k}(s_k),
%\label{Wvectm=1}
%\ee
%where we use the generalization of the prime circulator to
%vectorial index and argument.

Similarly to the Sylvester splitting theorem
(\ref{SylvWavesExpand}) we write the vector partition function
as a sum of {\it vector Sylvester waves}
\be
W({\bf s},{\bf D}^m) = \sum_{\bf j} W_{\bf j}({\bf s},{\bf D}^m) =
\sum_{j_1,\ldots,j_l} W_{\bf j}({\bf s},{\bf D}^m),
\label{SylvWavesExpandVect}
\ee
where ${\bf j}$ denotes $l$-dimensional vector $(j_1,\ldots,j_l)$.
The summation for each $j_k$ runs over all distinct factors of the
elements %of $i$-th row
of the matrix ${\bf D}^m$.
% The recursion (\ref{Wvectrecurs}) should hold for each
% "Sylvester vector wave" $W_{\bf k}({\bf s},{\bf D}^m)$.
%The polynomial part is made up of
%the waves for which the vectorial index ${\bf J}$ contains only zeros and
%ones. All other terms are periodic in one or more variables $s_j$.
% We assume that as in case of scalar restricted partition, the
%relations (\ref{Wvectrecurs}) and
% (\ref{Wvectsymm}) are valid for each vector wave.
% and, therefore, the polynomial part should satisfy them too.
The vector analog of the generator (\ref{generatorWj})
is written in the form (see \cite{BrionVergne1997})
$$
%\be
F_{\bf j}^m({\bf s},{\bf t}) =
\sum_{\rho_{j_1}, \rho_{j_2}, \ldots, \rho_{j_l}}
\frac{e^{\bf s \cdot t}\mbf{\rho}_{\bf j}^{-{\bf s}}}
{\prod_{i=1}^m (1-\mbf{\rho}_{\bf j}^{{\bf c}_i}e^{{-\bf c}_i \cdot {\bf t}})},
$$
%\label{generatorWvectk}
%\ee
where
\be
\mbf{\rho}_{\bf j}^{{\bf c}_i} = \prod_{k=1}^l \rho_{j_k}^{D_{ki}}.
\label{rhoelements}
\ee
It satisfies the relation similar to (\ref{Frecurs})
$$
%\be
F_{\bf j}^m({\bf s},{\bf t})-
F_{\bf j}^m({\bf s}-{\bf c}_m,{\bf t})=
F_{\bf j}^{m-1}({\bf s},{\bf t})
$$
%\label{Frecursvect}
%\ee
for any generator $F_{\bf j}^m({\bf s},{\bf t})$.
The following relation aslo holds
$$
F_{\bf j}^m({\bf s},{\bf t})=
(-1)^m 
F_{\bf j}^m(-{\bf s}-\mbf{\sigma}({\bf D}^m),-{\bf t}).
$$
The multidimensional residue of the generator
gives $W_{\bf j}({\bf s},{\bf D}^m)$ as the coefficient
of ${\bf t}^{-{\bf 1}} = \prod_{k=1}^l t_k^{-1}$, and
each vector wave satisfies the relations (\ref{Wvectrecurs}) and
(\ref{Wvectsymm}).

Let an equality $\mbf{\rho}_{\bf j}^{{\bf c}_i}=1$ holds for
$\omega_{\bf j}$ columns of the matrix ${\bf D}^m$, and
sort the matrix in such way that these
$\omega_{\bf j}$ columns come first, so that
$1 \le i \le \omega_{\bf j}$.
Introduce a homogeneous polynomial of degree $m$
$$
%\be
P_m({\bf t},{\bf D}^m)=\prod_{i=1}^m({\bf c}_i\cdot{\bf t})=
\sum_{|{\bf N}|=m} C_{\bf N}({\bf D}^m) {\bf t^N},
$$
%\label{factorP}
%\ee
and construct the modified generator
$$
\tilde{F}_{\bf j}^m({\bf s},{\bf t}) =
P_m({\bf t},{\bf D}^m)
F_{\bf j}^m({\bf s},{\bf t}),
$$
%in expansion of which we should select the coefficients
%of ${\bf t^n}$ with $|{\bf n}|=m-l$.
%
% The modified generator $\tilde{F}_{\bf k}^m({\bf s},{\bf t})$
which can be written as
\bea
%&&
\tilde{F}_{\bf j}^m({\bf s},{\bf t}) &=&
\sum_{\rho_{j_n}}
\mbf{\rho}_{\bf j}^{-{\bf s}}
\prod_{i=1}^{\omega_{\bf j}}
\frac{{\bf c}_i\cdot {\bf t}}{e^{{\bf c}_i\cdot {\bf t}}-1}
\times
e^{({\bf s}+ \mbf{\sigma}({\bf D}^m))\cdot {\bf t}}
\prod_{i=\omega_{\bf j}+1}^m
\frac{1-\mbf{\rho}_{\bf j}^{{\bf c}_j}}
{e^{{\bf c}_i \cdot {\bf t}}-\mbf{\rho}_{\bf j}^{{\bf c}_i}}
\cdot
\frac{{\bf c}_i\cdot {\bf t}}{1-\mbf{\rho}_{\bf j}^{{\bf c}_i}}
\nonumber \\
%\label{Wvectk1}  \\
&=&
\sum_{\rho_{j_n}}
\mbf{\rho}_{\bf j}^{-{\bf s}-\mbf{\sigma}({\bf D}^{m-\omega_{\bf j}})}
\sum_{{\bf n}_1}
B_{{\bf n}_1}^{(l,\omega_{\bf j})}(0|{\bf D}^{\omega_{\bf j}})
\frac{{\bf t}^{{\bf n}_1}}{{\bf n}_1!}
\nonumber \\
&\times&
\prod_{i=\omega_{\bf j}+1}^m
\frac{({\bf c}_i\cdot {\bf t})\mbf{\rho}_{\bf j}^{{\bf c}_i}}
{1-\mbf{\rho}_{\bf j}^{{\bf c}_i}}
\sum_{{\bf n}_2}
H_{{\bf n}_2}^{(l,m-\omega_{\bf j})}
({\bf s}+\mbf{\sigma}({\bf D}^m),
\mbf{\rho}_{\bf j}^{m-\omega_{\bf j}}|{\bf D}^{m-\omega_{\bf j}})
\frac{{\bf t}^{{\bf n}_2}}{{\bf n}_2!},
\nonumber
\eea
where elements of the vector $\mbf{\rho}^{m-\omega_{\bf j}}$
are $\mbf{\rho}_{\bf j}^{{\bf c}_i} \ne 1$ given by (\ref{rhoelements})
for $\omega_{\bf j}+1\le i\le m$.
Using the relation (\ref{BernEuler1vRel}) 
with $p_i$ such that
$\mbf{\rho}_{\bf j}^{p_i {\bf c}_i} = 1$
we can write
\bea
&&\prod_{i=\omega_{\bf j}+1}^m
\frac{({\bf c}_i\cdot {\bf t})\mbf{\rho}_{\bf j}^{{\bf c}_i}}
{1-\mbf{\rho}_{\bf j}^{{\bf c}_i}}
\sum_{{\bf n}_2}
H_{{\bf n}_2}^{(l,m-\omega_{\bf j})}
({\bf s}+\mbf{\sigma}({\bf D}^m),
\mbf{\rho}_{\bf j}^{{\bf c}_i}|{\bf D}^{m-\omega_{\bf j}})
\frac{{\bf t}^{{\bf n}_2}}{{\bf n}_2!}
\nonumber \\
&=&
\frac{1}{\pi({\bf p}^{m-\omega_{\bf j}})}
\sum_{{\bf n}_2}
\sum_{r_i=0}^{p_i-1}
\mbf{\rho}_{\bf j}^{-r_i{\bf c}_i}
B_{{\bf n}_2}^{(l,m-\omega_{\bf j})}
({\bf s}+%{\bf r \cdot \bf D}^{m-\omega_{\bf j}}|
\sum_{i=\omega_{\bf j}+1}^m r_i{\bf c}_i|
\{p_{\omega_{\bf j}+1} {\bf c}_{\omega_{\bf j}+1},\ldots,p_m {\bf c}_m\})
\frac{{\bf t}^{{\bf n}_2}}{{\bf n}_2!}.
\nonumber
\eea
Employing the binomial formula (\ref{Bern1vbinom})
we obtain $\tilde{F}_{\bf j}^m({\bf s},{\bf t})$
in the form
$$
%\be
\tilde{F}_{\bf j}^m({\bf s},{\bf t})=
\frac{1}{\pi({\bf p}^{m-\omega_{\bf j}})}
\sum_{{\bf n}}\sum_{r_i=0}^{p_i-1}
B_{\bf n}^{(l,m)}
({\bf s}+\mbf{\sigma}({\bf D}^m)+
{\bf r \cdot \bf D}^{m-\omega_{\bf j}}
%\!\!\!\sum_{i=\omega_{\bf j}+1}^m \!\!\!r_i{\bf c}_i
|{\bf D}^m_{\bf j})
\frac{{\bf t}^{{\bf n}}}{{\bf n}!}\cdot
\sum_{\rho_{j_i}}
\mbf{\rho}_{\bf j}^{-{\bf s}-
{\bf r \cdot \bf D}^{m-\omega_{\bf j}}
%\sum_{i=\omega_{\bf j}+1}^m (r_i+1){\bf c}_i
},
$$
%\label{Wvectk2}
%\ee
where ${\bf D}^m_{\bf j}$ is the ${\bf j}$-modified
matrix of the form
$$
{\bf D}^m_{\bf j} =\{{\bf c}_1, \ldots, {\bf c}_{\omega},
p_{\omega_{\bf j}+1}{\bf c}_{\omega_{\bf j}+1},\ldots,p_m{\bf c}_m\}.
$$
Noting that for the ${\bf j}$-modified matrix ${\bf D}^m_{\bf j}$
$$
P_m({\bf t},{\bf D}^m_{\bf j})=
\pi({\bf p}^{m-\omega_{\bf j}}) P_m({\bf t},{\bf D}^m), \ \ \
\mbf{\sigma}({\bf D}^m_{\bf j}) =
\mbf{\sigma}({\bf D}^m) +
\sum_{i=\omega_{\bf j}+1}^m (p_i-1) {\bf c}_i,
$$
and using the vector prime circulator notation 
$$
%\be
\Psi_{\bf c}({\bf s})=\prod_{k=1}^l \Psi_{c_k}(s_k),
$$
%\label{Psi_vector_def}
%\ee
we write for the
generator ${F}_{\bf j}^m({\bf s},{\bf t})$
$$
%\be
{F}_{\bf j}^m({\bf s},{\bf t})=
P_m^{-1}({\bf t},{\bf D}^m_{\bf j})
\sum_{{\bf n}}\sum_{r_i=0}^{p_i-1}
B_{\bf n}^{(l,m)}
({\bf s}+\mbf{\sigma}({\bf D}^m_{\bf j})-
{\bf r \cdot \bf D}^{m-\omega_{\bf j}}|
%\!\!\!\sum_{j=\omega_{\bf j}+1}^m \!\!\! r_i{\bf c}_i|
{\bf D}^m_{\bf j})
\Psi_{\bf j}
({\bf s}-
{\bf r \cdot \bf D}^{m-\omega_{\bf j}})
%\!\!\!\sum_{i=\omega_{\bf j}+1}^m \!\!\! r_i{\bf c}_i)
\frac{{\bf t}^{{\bf n}}}{{\bf n}!}.
$$
%\label{Wvectk2a}
%\ee
In the above expression
only terms with $|{\bf n}|=m-l$ contribute to the vector Sylvester wave
$W_{\bf j}({\bf s},{\bf D}^m)$, which is found also as a constant term of
${F}_{\bf j}^m({\bf s},{\bf t}){\bf t}^{\bf 1}$, equal to a fraction with both
numerator and denominator being homogeneous polynomials of degree $m$
$$
%\be
W_{\bf j}({\bf s},{\bf D}^m) =
\lim_{{\bf t}\to{\bf 0}}
P_m^{-1}({\bf t},{\bf D}^m_{\bf j})\!\!
\sum_{|{\bf n}|=m-l}\sum_{r_i=0}^{p_i-1}
B_{{\bf n}}^{(l,m)}
%({\bf s}+\mbf{\sigma}({\bf D}^m_{\bf j})-
%{\bf r \cdot \bf D}^{m-\omega_{\bf j}}|{\bf D}^m_{\bf j})
({\bf s}'|{\bf D}^m_{\bf j})
\Psi_{\bf j}({\bf s}')
%\Psi_{\bf j}({\bf s}+\mbf{\sigma}({\bf D}^m_{\bf j})-
%{\bf r \cdot \bf D}^{m-\omega_{\bf j}})
\frac{{\bf t}^{{\bf n}+{\bf 1}}}{{\bf n}!},
\ \
%$$
%\label{Wvectk2b}
%\ee
%where
%$$
{\bf s}'={\bf s}+\mbf{\sigma}({\bf D}^m_{\bf j})-
{\bf r \cdot \bf D}^{m-\omega_{\bf j}}.
%\!\!\!\sum_{i=\omega_{\bf j}+1}^m \!\!\! r_i{\bf c}_i.
$$

%It should be noted that each vector ${\bf n}$ with $|{\bf n}|=m-l$
%in (\ref{Wvectk2}) corresponds to separate chamber in the chamber complex
%of quasipolynomiality. Thus 
It is convenient to write each vector
Sylvester wave $W_{\bf j}({\bf s},{\bf D}^m)$ as a sum of %corresponding
quasipolynomials, which we call {\it partial vector
Sylvester waves}:
\be
W_{\bf j}({\bf s},{\bf D}^m) =
%\cup_{|{\bf n}|=m-l}
\sum_{|{\bf n}|=m-l}
W_{\bf j}^{\bf n}({\bf s},{\bf D}^m).
\label{Wvectk_union}
\ee
The number of different partial waves is found as a number of 
representations of an integer $m-l$ as a sum of
$l$ nonnegative integers, i.e., it is equal to
$$
W(m-l,{\bf 1}^l)=\binom{m-1}{m-l}.
$$
The partial wave
$W_{\bf j}^{\bf n}({\bf s},{\bf D}^m)$ is found as
($|{\bf n}|=m-l$)
$$
%\be
W_{\bf j}^{\bf n}({\bf s},{\bf D}^m) =
\lim_{{\bf t}\to{\bf 0}}
P_m^{-1}({\bf t},{\bf D}^m_{\bf j})
\sum_{r_i=0}^{p_i-1}
B_{{\bf n}}^{(l,m)}
({\bf s}+\mbf{\sigma}({\bf D}^m_{\bf j})-
{\bf r \cdot \bf D}^{m-\omega_{\bf j}}
%\!\!\!\sum_{i=\omega_{\bf j}+1}^m \!\!\! r_i{\bf c}_i
|{\bf D}^m_{\bf j})
\Psi_{\bf j}({\bf s}-
{\bf r \cdot \bf D}^{m-\omega_{\bf j}})
%\!\!\!\sum_{i=\omega_{\bf j}+1}^m \!\!\! r_i{\bf c}_i)
\frac{{\bf t}^{{\bf n}+{\bf 1}}}{{\bf n}!}.
$$
%\label{Wvectpartial}
%\ee
Using a parametrization ${\bf t}=\mbf{\alpha}t$ we compute the above limit
as
\bea
W_{\bf j}^{\bf n}({\bf s},{\bf D}^m)  &=&
P_m^{-1}(\mbf{\alpha},{\bf D}^m_{\bf j})
\sum_{r_i=0}^{p_i-1}
B_{{\bf n}}^{(l,m)}
({\bf s}+\mbf{\sigma}({\bf D}^m_{\bf j})-
{\bf r \cdot \bf D}^{m-\omega_{\bf j}}
%\!\!\!\sum_{i=\omega_{\bf j}+1}^m \!\!\! r_i{\bf c}_i
|{\bf D}^m_{\bf j})
\Psi_{\bf j}({\bf s}-
{\bf r \cdot \bf D}^{m-\omega_{\bf j}})
%\!\!\!\sum_{i=\omega_{\bf j}+1}^m \!\!\! r_i{\bf c}_i)
\frac{\mbf{\alpha}^{{\bf n}+{\bf 1}}}{{\bf n}!}%.
\nonumber \\
&=&
C_{\bf n}(\mbf{\alpha},{\bf D}^m_{\bf j})
\sum_{r_i=0}^{p_i-1}
B_{{\bf n}}^{(l,m)}
({\bf s}+\mbf{\sigma}({\bf D}^m_{\bf j})-
{\bf r \cdot \bf D}^{m-\omega_{\bf j}}
%\!\!\!\sum_{i=\omega_{\bf j}+1}^m \!\!\! r_i{\bf c}_i
|{\bf D}^m_{\bf j})
\Psi_{\bf j}({\bf s}-
{\bf r \cdot \bf D}^{m-\omega_{\bf j}}),
%\!\!\!\sum_{i=\omega_{\bf j}+1}^m \!\!\! r_i{\bf c}_i),
%\label{Wvectpartial1}
\nonumber
\eea
where 
$$
%\be
C_{\bf n}(\mbf{\alpha},{\bf D}^m_{\bf j}) =
\frac{\mbf{\alpha}^{{\bf n}+{\bf 1}}}
{{\bf n}!P_m(\mbf{\alpha},{\bf D}^m_{\bf j})},
\ \ \ |{\bf n}|=m-l.
$$
%\label{polynomialC}
%\ee

We define the polynomial part of the vector partition
function for ${\bf j}={\bf 1}$, the corresponding
partial polynomial is equal to
\be
W_{\bf 1}^{\bf n}({\bf s},{\bf D}^m) =
C_{\bf n}(\mbf{\alpha},{\bf D}^m)
B_{{\bf n}}^{(l,m)}
({\bf s}+\mbf{\sigma}({\bf D}^m)|{\bf D}^m).
%\frac{\mbf{\alpha}^{{\bf n}+{\bf 1}}}{{\bf n}!}
\label{Wvect1}
\ee
The partial polynomial part for
the ${\bf j}$-modified matrix ${\bf D}^m_{\bf j}$ reads
$$
%\be
W_{\bf 1}^{\bf n}({\bf s},{\bf D}^m_{\bf j}) =
C_{\bf n}(\mbf{\alpha},{\bf D}^m_{\bf j})
B_{{\bf n}}^{(l,m)}
({\bf s}+\mbf{\sigma}({\bf D}^m_{\bf j})|{\bf D}^m_{\bf j})
%\frac{\mbf{\alpha}^{{\bf n}+{\bf 1}}}{{\bf n}!}
,
$$
%\label{Wvect1k}
%\ee
and similarly to the scalar case the partial
vector Sylvester wave $W_{\bf j}^{\bf n}({\bf s},{\bf D}^m)$
for arbitrary ${\bf j} \ne {\bf 1}$ can be written as
a linear superposition
of the partial polynomial part for the % ${\bf k}$-modified
matrix ${\bf D}^m_{\bf j}$ multiplied by the
corresponding prime circulator:
\be
W_{\bf j}^{\bf n}({\bf s},{\bf D}^m) =
\sum_{r_i=0}^{p_i-1}
W_{\bf 1}^{\bf n}({\bf s}-
{\bf r \cdot \bf D}^{m-\omega_{\bf j}}
%\!\sum_{i=\omega_{\bf j}+1}^m \!\!\!r_i{\bf c}_i
,{\bf D}^m_{\bf j})
\Psi_{\bf j}
({\bf s}-
{\bf r \cdot \bf D}^{m-\omega_{\bf j}}
%\!\sum_{i=\omega_{\bf j}+1}^m \!\!\!r_i{\bf c}_i
).
\label{Wvectkthrough1}
\ee

Combining (\ref{SylvWavesExpandVect},\ref{Wvectk_union}) and 
(\ref{Wvectkthrough1}) we arrive at the final expression of
the restricted vector partition function as a linear superposition of the 
vector Bernoulli polynomials of higher order multiplied by the 
vector prime circulators:
\be
W({\bf s},{\bf D}^m)=
\sum_{\bf j}
\sum_{|{\bf n}|=m-l}
\sum_{r_i=0}^{p_i-1}
W_{\bf 1}^{\bf n}({\bf s}-
{\bf r \cdot \bf D}^{m-\omega_{\bf j}}
%\!\!\!\sum_{i=\omega_{\bf j}+1}^m \!\!\!r_i{\bf c}_i
,{\bf D}^m_{\bf j})
\Psi_{\bf j}
({\bf s}-
{\bf r \cdot \bf D}^{m-\omega_{\bf j}}
%\!\!\!\sum_{i=\omega_{\bf j}+1}^m \!\!\!r_i{\bf c}_i
),
\label{Wvectfinal}
\ee
where
$W_{\bf 1}^{\bf n}({\bf s},{\bf D}^m_{\bf j})$ is given by
\be
W_{\bf 1}^{\bf n}({\bf s},{\bf D}^m_{\bf j})=
\frac{\mbf{\alpha}^{{\bf n}+{\bf 1}}}
{{\bf n}!P_m(\mbf{\alpha},{\bf D}^m_{\bf j})}\;
B_{{\bf n}}^{(l,m)}
({\bf s}+\mbf{\sigma}({\bf D}^m_{\bf j})|{\bf D}^m_{\bf j}).
\label{Wvect1kC}
\ee
It also can be written in the form
of a mixture of partial waves
%$$
\be
W({\bf s},{\bf D}^m)=\!\!
\sum_{|{\bf n}|=m-l}\!\!
C_{\bf n}(\mbf{\alpha},{\bf D}^m)
W^{\bf n}({\bf s},{\bf D}^m),
\label{Wvectfinalpartial}
%$$
\ee
where the superposition coefficients depend on the vector $\mbf{\alpha}$
and the partial wave is given by
$$
%\be
W^{\bf n}({\bf s},{\bf D}^m)=
\sum_{\bf j}
\sum_{r_i=0}^{p_i-1}
\pi^{-1}({\bf p}^{m-\omega_{\bf j}})
B_{{\bf n}}^{(l,m)}
({\bf s}+\mbf{\sigma}({\bf D}^m_{\bf j})-
{\bf r \cdot \bf D}^{m-\omega_{\bf j}}
|{\bf D}^m_{\bf j})
\Psi_{\bf j}({\bf s}-
{\bf r \cdot \bf D}^{m-\omega_{\bf j}}).
$$
%\label{Wvectpartial_def}
%\ee

In the scalar case we have $l=1$, ${\bf D}^m={\bf d}^m, n=m-1$,
and
$
P_m(\alpha,{\bf d}^m_{j})=\pi({\bf d}^m_{j})\alpha^m.
$
Thus, the expression (\ref{Wvect1kC}) reduces to (\ref{W1j}), 
and noting that in
the scalar case all $p_i=j$, one finds that (\ref{Wvectkthrough1}) transforms
into
(\ref{WjthroughW1j}).

We present several examples of application of the formula
(\ref{Wvectfinal}) and show that it gives not only the vector partition
function but also enables to find the chamber structure of the system.

%%%%%%%%%%%%%%%%%%%%%%%%%%%%%%%%%%%%%%%%%%%%%%%%%%%%%%%%%%%%%%%%%%%%%%%%%%%%%%%%%
\subsection{Example 1}
Consider computation of the vector partition function
for the matrix $(m=3,l=2)$
$$
{\bf D}^3 =
\left(
\begin{array}{ll}
1 \ 2 \ 0 \\
1 \ 0 \ 1
\end{array}
\right), \ \ \ \
{\bf s} = (s_1,s_2),
$$
first using the partial fraction expansion suggested in \cite{Beck2004}.
An idea of the method is based on the definition of the generating function
(\ref{WvectGF}). It is easily seen that
$W({\bf s},{\bf D})$ can be computed as the constant term in the expansion
of the following expression:
\be
W({\bf s},{\bf D}^m)=\mbox{const}\;
\left[\frac{1}{{\bf t^s}} \prod_{i=1}^m 
\frac{1}{1-{\bf t}^{{\bf c}_i}} \right].
\label{Beck1}
\ee
Assuming all but one (say $t_1$) components of the vector ${\bf t}$ to be
constant expand the r.h.s. in (\ref{Beck1}) into the partial fractions in $t_1$.
This expansion contains both analytic and meromorphic parts w.r.t $t_1=0$.
The meromorphic part doesn't contribute to the $t_1$-constant term, so it can be
dropped. The constant term of the analytic part depends on the remaining 
components of ${\bf t}$, so that such expansion eliminates $t_1$. Applying this
procedure repeatedly we eliminate all components of ${\bf t}$ and
arrive to the result.

%{\bf Case 1.}
Assume $t_1$ constant and make partial fraction expansion w.r.t.
$t_2$, obtaining
$$
\frac{1}{(1-t_2)(1-t_1t_2)t_2^{s_2}} =
\frac{1}{1-t_1}\left[\frac{1}{1-t_2} - \frac{t_1^{s_2+1}}{1-t_1t_2} \right] +
\mbox{MMP},
$$
where MMP stands for the meromorphic part.
Thus we have
\bea
W({\bf s},{\bf D}^3) &=&
\mbox{const}_{t_1} \left[ \frac{1}{(1-t_1^2)(1-t_1)t_1^{s_1}}
\mbox{const}_{t_2} \left[\frac{1}{1-t_2} - \frac{t_1^{s_2+1}}{1-t_1t_2} \right]
\right]
\nonumber \\
&=&
\mbox{const}_{t} \left[
\frac{1-t^{s_2+1}}{(1-t^2)(1-t)t^{s_1}}\right].
\nonumber
%\label{Way1a}
\eea
Finding that for $a \ge 0$
$$
\mbox{const}_{t} \frac{1}{(1-t^2)(1-t)t^{a}} =
\mbox{const}_{t} \left[
\frac{1}{2(1-t)^2}+\frac{1+2a}{4(1-t)}+\frac{(-1)^a}{4(1+t)} + \mbox{MMP}
\right] =
\frac{a}{2}+\frac{3+(-1)^a}{4},
$$
we obtain (note that $s_1 \ge 0$)
$$
\mbox{const}_{t} \left[
\frac{1}{(1-t^2)(1-t)t^{s_1}}\right] = \frac{s_1}{2}+\frac{3+(-1)^{s_1}}{4},
$$
and
$$
\mbox{const}_{t} \left[
\frac{1}{(1-t^2)(1-t)t^{s_1-s_2-1}}\right] =
\left\{
\begin{array}{ll}
0,                                               & s_1-s_2-1 < 0, \\
\frac{s_1-s_2}{2} + \frac{1-(-1)^{s_1-s_2}}{4},  & s_1-s_2-1 \ge 0.
\end{array}
\right.
$$
Combining these results we arrive at the final expression for the
partition function
\be
W({\bf s},{\bf D}^3) =
\left\{
\begin{array}{ll}
\frac{s_1}{2}+\frac{3+(-1)^{s_1}}{4}, & s_1-s_2-1 < 0, \\
\frac{s_2+1}{2}+\frac{(-1)^{s_1}+(-1)^{s_1-s_2}}{4}, & s_1-s_2-1 \ge 0.
\end{array}
\right.
\label{Way1res}
\ee
% The boundary between the two chambers is the line $s_1=s_2$.
It is easy to
check that the expressions in (\ref{Way1res}) coincide
for $s_1=s_2$ and $s_1=s_2-1$.

Turning to the formula (\ref{Wvectfinalpartial}) and having
$|{\bf n}|=1$, we find that there are two
partial waves corresponding to
$
{\bf n}_1=(1,0),
{\bf n}_2=(0,1).
$
Computing the polynomial
$P_3(\mbf{\alpha},{\bf D}^3)=
2\alpha_1\alpha_2(\alpha_1+\alpha_2)$
we obtain the coefficients:
$$
C_{(1,0)}(\mbf{\alpha},{\bf D}^3) =
\frac{\alpha_1}{2(\alpha_1+\alpha_2)}, \ \ \
C_{(0,1)}(\mbf{\alpha},{\bf D}^3) =
\frac{\alpha_2}{2(\alpha_1+\alpha_2)}.
$$
Finding $\mbf{\sigma}({\bf D}^3)=(3,2)$
and using (\ref{Wvect1}) we
have for the partial polynomial parts
\be
W_{\bf 1}^{(1,0)}({\bf s},{\bf D}^3) =
C_{(1,0)}(\mbf{\alpha},{\bf D}^3)
\left(s_1+\frac{3}{2}\right),\ \
%\nonumber \\
W_{\bf 1}^{(0,1)}({\bf s},{\bf D}^3) =
C_{(0,1)}(\mbf{\alpha},{\bf D}^3)
(s_2 + 1).
\label{Wvect1example0}
\ee
In addition to the polynomial contributions there
are nonzero terms
corresponding to ${\bf j}_1=(2,1)$
and ${\bf j}_2=(2,2)$, while the
term with ${\bf j}=(1,2)$
doesn't contribute to the final result.

The ${\bf j}_1$-modified matrix ${\bf D}^m_{{\bf j}_1}$ reads
($p_3=2$)
$$
{\bf D}^3_{(2,1)} =
\left(
\begin{array}{ll}
2 \ 0 \ 2  \\
0 \ 1 \ 2
\end{array}
\right),
$$
where the columns are sorted already.
In order to use % the formula
(\ref{Wvectkthrough1}) we first compute the partial polynomials for the
${\bf j}_1$-modified matrix:
$$
%\be
W_{\bf 1}^{(1,0)}({\bf s},{\bf D}^3_{(2,1)}) =
C_{(1,0)}(\mbf{\alpha},{\bf D}^3)
\left(\frac{s_1}{2}+1\right), \ \
%\nonumber \\
W_{\bf 1}^{(0,1)}({\bf s},{\bf D}^3_{(2,1)}) =
C_{(0,1)}(\mbf{\alpha},{\bf D}^3)
\left(\frac{s_2}{2} + \frac{3}{4} \right),
$$
%\label{Wvect1kexample0}
%\ee
and find
$$
W_{(2,1)}^{{\bf n}_i}({\bf s},{\bf D}^3) =
\sum_{r=0}^{1}
W_{\bf 1}^{{\bf n}_i}((s_1-r,s_2-r),{\bf D}^3_{(2,1)})
\Psi_2(s_1-r),
$$
arriving at
\be
W_{(2,1)}^{(1,0)}({\bf s},{\bf D}^3) =
\frac{1}{2}C_{(1,0)}(\mbf{\alpha},{\bf D}^3)
\Psi_2(s_1),\ \
%\nonumber \\
W_{(2,1)}^{(0,1)}({\bf s},{\bf D}^3) =
\frac{1}{2}C_{(0,1)}(\mbf{\alpha},{\bf D}^3)
\Psi_2(s_1).
\label{Wvectkkexample0}
\ee
The ${\bf j}_2$-modified matrix ${\bf D}^m_{{\bf j}_2}$ reads
($p_2=p_3=2$)
$$
{\bf D}^3_{(2,2)} =
\left(
\begin{array}{ll}
2 \ 0 \ 2  \\
0 \ 2 \ 2
\end{array}
\right).
$$
%where the columns are sorted already.
%In order to use % the formula
%(\ref{Wvectkthrough1}) we first compute 
Repeating the computation we start with
the partial polynomials for the
${\bf j}_2$-modified matrix:
$$
%\be
W_{\bf 1}^{(1,0)}({\bf s},{\bf D}^3_{(2,2)}) =
C_{(1,0)}(\mbf{\alpha},{\bf D}^3)
\left(\frac{s_1}{4}+\frac{1}{2}\right),\ \
%\nonumber \\
W_{\bf 1}^{(0,1)}({\bf s},{\bf D}^3_{(2,2)}) =
C_{(0,1)}(\mbf{\alpha},{\bf D}^3)
\left(\frac{s_2}{4} + \frac{1}{2} \right).
$$
%\label{Wvect1kexample0a}
%\ee
The formula (\ref{Wvectkthrough1}) takes form
$$
W_{(2,2)}^{{\bf n}_i}({\bf s},{\bf D}^3) =
\sum_{r_1,r_2=0}^{1}
W_{\bf 1}^{{\bf n}_i}((s_1-r_1,s_2-r_1-r_2),{\bf D}^3_{(2,2)})
\Psi_2(s_1-r_1)\Psi_2(s_2-r_1-r_2),
$$
and we obtain
\be
W_{(2,2)}^{(1,0)}({\bf s},{\bf D}^3) =
0,\ \
%\nonumber \\
W_{(2,2)}^{(0,1)}({\bf s},{\bf D}^3) =
\frac{1}{2}C_{(0,1)}(\mbf{\alpha},{\bf D}^3)
\Psi_2(s_1)\Psi_2(s_2).
\label{Wvectkkexample0a}
\ee
Combining the expressions (\ref{Wvect1example0}),
(\ref{Wvectkkexample0}) and
(\ref{Wvectkkexample0a}) we arrive at the final result
\be
W({\bf s},{\bf D}^3)=
%\frac{\alpha_1}{2(\alpha_1+\alpha_2)}
C_{(1,0)}(\mbf{\alpha},{\bf D}^3)
W^{(1,0)}({\bf s},{\bf D}^3) +
%\frac{\alpha_2}{2(\alpha_1+\alpha_2)}
C_{(0,1)}(\mbf{\alpha},{\bf D}^3)
W^{(0,1)}({\bf s},{\bf D}^3),
\label{example0}
\ee
where
\bea
W^{(1,0)}({\bf s},{\bf D}^3) &=&
s_1+\frac{3}{2}+
\frac{1}{2}\Psi_2(s_1),
\nonumber \\
W^{(0,1)}({\bf s},{\bf D}^3) &=&
s_2+ 1 +
\frac{1}{2}\Psi_2(s_1)+
\frac{1}{2}\Psi_2(s_1-s_2).
%\label{Wvectfinalexample0}
\nonumber
\eea
Noting that $\Psi_2(s)=(-1)^s$ one can check by straightforward computation
that the case $\alpha_2/\alpha_1=0$ produces the first line in
(\ref{Way1res}), while $\alpha_1/\alpha_2=0$ gives the last line.
We see that
proper choice of the vector $\mbf{\alpha}$ corresponds to the selection
of the chamber.

%%%%%%%%%%%%%%%%%%%%%%%%%%%%%%%%%%%%%%%%%%%%%%%%%%%%%%%%%%%%%%%%%%%%%%%%%%%%%%%%%
\subsection{Example 2}
Consider another example for $l=2, m=4$ and
$$
{\bf D}^4 =
\left(
\begin{array}{ll}
1 \ 2 \ 1 \ 0 \\
1 \ 1 \ 0 \ 1
\end{array}
\right), \ \ \ \
{\bf s} = (s_1,s_2).
$$
The straightforward computation using the partial fraction expansion gives the 
following result (see \cite{Beck2004}):
\be
W({\bf s},{\bf D}^4) =
\left\{
\begin{array}{ll}
\frac{s_1^2}{4}+s_1+\frac{7+(-1)^{s_1}}{8}, &  s_1 \le s_2,
\nonumber \\
s_1s_2-\frac{s_1^2+2s_2^2}{4}+\frac{s_1+s_2}{2}+\frac{7+(-1)^{s_1}}{8}, & 
s_1/2-1 \le s_2 \le s_1+1, \\
\frac{s_2^2}{2}+\frac{3s_2}{2}+1, & s_2 \le s_1/2.
\end{array}
\right.
\label{Beck_ex_res}
\ee
As $|{\bf n}|=2$, we have three
partial waves corresponding to
$
{\bf n}_1=(2,0),
{\bf n}_2=(1,1),
{\bf n}_3=(0,2).
$
Computing the polynomial
$P_4(\mbf{\alpha},{\bf D}^4)=
\alpha_1\alpha_2(\alpha_1+\alpha_2)(2\alpha_1+\alpha_2)$
we find the coefficients:
$$
C_{(2,0)}(\mbf{\alpha},{\bf D}^4) = %\frac{1}{4}
\frac{\alpha_1^2}{2(\alpha_1+\alpha_2)(2\alpha_1+\alpha_2)},
$$
$$
C_{(1,1)}(\mbf{\alpha},{\bf D}^4) = %\frac{1}{3}
\frac{\alpha_1\alpha_2}{(\alpha_1+\alpha_2)(2\alpha_1+\alpha_2)},
$$
$$
C_{(0,2)}(\mbf{\alpha},{\bf D}^4) = %\frac{1}{2}
\frac{\alpha_2^2}{2(\alpha_1+\alpha_2)(2\alpha_1+\alpha_2)}.
$$
Finding $\mbf{\sigma}({\bf D}^4)=(4,3)$ 
and using (\ref{Wvect1}) we
have for the partial polynomial parts
\bea
W_{\bf 1}^{(2,0)}({\bf s},{\bf D}^4) &=&
%\frac{1}{4}
%C_{(3,1)}(\mbf{\alpha},{\bf D}^4)
%B_{(2,0)}^{(2,4)}((s_1+4,s_2+3)|{\bf D}^4) =
C_{(2,0)}(\mbf{\alpha},{\bf D}^4)
\left(s_1^2+4s_1+\frac{7}{2}\right),
%\left(\frac{s_1^2}{4}+s_1+\frac{7}{8}\right),
\nonumber \\
W_{\bf 1}^{(1,1)}({\bf s},{\bf D}^4) &=&
%\frac{1}{3}
%C_{(2,2)}(\mbf{\alpha},{\bf D}^4)
%B_{(1,1)}^{(2,4)}((s_1+4,s_2+3)|{\bf D}^4) =
C_{(1,1)}(\mbf{\alpha},{\bf D}^4)
\left(s_1s_2 +
\frac{3s_1}{2} + 2 s_2 + \frac{11}{4}\right),
%\left(\frac{s_1s_2}{3} +
%\frac{s_1}{2} + \frac{2 s_2}{3} + \frac{11}{12}\right),
\label{Wvect1example} \\
W_{\bf 1}^{(0,2)}({\bf s},{\bf D}^4) &=&
% \frac{1}{2}
%C_{(1,3)}(\mbf{\alpha},{\bf D}^4)
%B_{(0,2)}^{(2,4)}((s_1+4,s_2+3)|{\bf D}^4) =
C_{(0,2)}(\mbf{\alpha},{\bf D}^4)
\left(s_2^2+3s_2+2\right).
%\left(\frac{s_2^2}{2}+\frac{3s_2}{2}+1\right).
\nonumber
\eea
In addition to the polynomial part only the
term corresponding to ${\bf j}=(2,1)$ produces nonzero contribution to the
vector partition function, while two other terms with
${\bf j}=(1,2)$ and ${\bf j}=(2,2)$ don't contribute into the result.
The ${\bf j}$-modified matrix ${\bf D}^m_{\bf j}$ reads
($p_3=p_4=2$)
$$
{\bf D}^4_{(2,1)} =
\left(
\begin{array}{ll}
2 \ 0 \ 2 \ 2 \\
1 \ 1 \ 2 \ 0
\end{array}
\right).
$$
% where the columns are sorted already.
We first compute the partial polynomials
for the
${\bf j}$-modified matrix:
\bea
W_{\bf 1}^{(2,0)}({\bf s},{\bf D}^4_{(2,1)}) &=&
%\frac{1}{16}
%B_{(2,0)}^{(2,4)}
%((s_1+6,s_2+4)|{\bf D}^4_{(2,1)}) =
C_{(2,0)}(\mbf{\alpha},{\bf D}^4)
\left(\frac{s_1^2}{4}+\frac{3s_1}{2}+ 2\right),
% \frac{s_1^2}{16}+\frac{3s_1}{8}+ \frac{1}{2},
\nonumber \\
W_{\bf 1}^{(1,1)}({\bf s},{\bf D}^4_{(2,1)}) &=&
%\frac{1}{12}
%B_{(1,1)}^{(2,4)}((s_1+6,s_2+4)|{\bf D}^4_{(2,1)}) =
C_{(1,1)}(\mbf{\alpha},{\bf D}^4)
\left(\frac{s_1s_2}{4} + \frac{s_1}{2} + \frac{3s_2}{4} + \frac{11}{8}\right),
%\frac{s_1s_2}{12} + \frac{s_1}{6} + \frac{s_2}{4} + \frac{11}{24},
\nonumber \\
%\label{Wvect1kexample} \\
W_{\bf 1}^{(0,2)}({\bf s},{\bf D}^4_{(2,1)}) &=&
%\frac{1}{8}
%B_{(0,2)}^{(2,4)}
%((s_1+6,s_2+4)|{\bf D}^4_{(2,1)}) =
C_{(0,2)}(\mbf{\alpha},{\bf D}^4)
\left(\frac{s_2^2}{4}+s_2+\frac{7}{8}\right),
%\frac{s_2^2}{8}+\frac{s_2}{2}+\frac{7}{16}.
\nonumber
\eea
and use (\ref{Wvectkthrough1}) to have
$$
W_{(2,1)}^{{\bf n}_i}({\bf s},{\bf D}^4) =
\sum_{r_3=0}^{1}\sum_{r_4=0}^{1}
W_{\bf 1}^{{\bf n}_i}((s_1-r_3-r_4,s_2-r_3),{\bf D}^4_{(2,1)})
\Psi_2(s_1-r_3-r_4),
$$
and obtain
\bea
W_{(2,1)}^{(2,0)}({\bf s},{\bf D}^4) &=&
\frac{1}{2}C_{(2,0)}(\mbf{\alpha},{\bf D}^4)
\Psi_2(s_1),
%\frac{1}{8}\Psi_2(s_1),
\nonumber \\
W_{(2,1)}^{(1,1)}({\bf s},{\bf D}^4) &=&
\frac{1}{4}C_{(1,1)}(\mbf{\alpha},{\bf D}^4)\Psi_2(s_1),
%\frac{1}{12}\Psi_2(s_1),
\label{Wvectkkexample} \\
W_{(2,1)}^{(0,2)}({\bf s},{\bf D}^4) &=&
0.
\nonumber
\eea
Combining the expressions (\ref{Wvect1example}) and
(\ref{Wvectkkexample}) we arrive at the final result
$$
%\be
W({\bf s},{\bf D}^4)=
C_{(2,0)}(\mbf{\alpha},{\bf D}^4)
W^{(2,0)}({\bf s},{\bf D}^4) +
C_{(1,1)}(\mbf{\alpha},{\bf D}^4)
W^{(1,1)}({\bf s},{\bf D}^4) +
C_{(0,2)}(\mbf{\alpha},{\bf D}^4)
W^{(0,2)}({\bf s},{\bf D}^4),
$$
%\label{Wvectfinalres}
%\ee
with
\bea
W^{(2,0)}({\bf s},{\bf D}^4) &=&
%C_{(2,0)}(\mbf{\alpha},{\bf D}^4)
%\frac{\alpha_1^2}{2(\alpha_1+\alpha_2)(2\alpha_1+\alpha_2)}
%\left(
s_1^2+4s_1+\frac{7}{2}+
\frac{1}{2}\Psi_2(s_1),
%\right),
%\frac{s_1^2}{4}+s_1+\frac{7}{8}+\frac{1}{8}\Psi_2(s_1),
\nonumber \\
W^{(1,1)}({\bf s},{\bf D}^4) &=&
%C_{(1,1)}(\mbf{\alpha},{\bf D}^4)
%\frac{\alpha_1\alpha_2}{(\alpha_1+\alpha_2)(2\alpha_1+\alpha_2)}
%\left(
s_1s_2 +
\frac{3s_1}{2} + 2 s_2 + \frac{11}{4}+
\frac{1}{4}\Psi_2(s_1),
%\right),
%\frac{s_1s_2}{3} + \frac{s_1}{2} + \frac{2 s_2}{3} +
%\frac{11}{12}+\frac{1}{12}\Psi_2(s_1),
\nonumber \\
%\label{Wvectfinalexample} \\
W^{(0,2)}({\bf s},{\bf D}^4) &=&
%C_{(0,2)}(\mbf{\alpha},{\bf D}^4)
%\frac{\alpha_2^2}{2(\alpha_1+\alpha_2)(2\alpha_1+\alpha_2)}
%\left(
s_2^2+3s_2+2.
%\right).
%\frac{s_2^2}{2}+\frac{3s_2}{2}+1.
\nonumber
\eea
%Noting that $\Psi_2(s)=(-1)^s$ one can check by straightforward computation
We find
that the case $\alpha_2/\alpha_1=0$ produces the first line in
(\ref{Beck_ex_res}), while $\alpha_1/\alpha_2=0$ corresponds to the last line.
The second line is obtained as a 
\underline{real} part of $W({\bf s},{\bf D}^4)$
for $\alpha_2/\alpha_1=-1 \pm i$.
Thus, we see again that
the choice of the vector $\mbf{\alpha}$ corresponds to
selection of the chamber.

%%%%%%%%%%%%%%%%%%%%%%%%%%%%%%%%%%%%%%%%%%%%%%%%%%%%%%%%%%%%%%%%%%%%%%%%%%%%%%%%%
\subsection{Example 3}
It should be noted that the formula (\ref{Wvectfinal})
remains valid also in the %degenerate
case $m=l$. To illustrate it consider the case $m=l=2$ with
$$
{\bf D}^2 =
\left(
\begin{array}{ll}
1 \ 2  \\
1 \ 0
\end{array}
\right), \ \ \ \
{\bf s} = (s_1,s_2).
$$
The solution of the corresponding linear problem reads
$x_1=s_2,\ x_2=(s_1-s_2)/2$, thus the number of solutions equals to
one for even nonnegative differences $s_1-s_2$ and zero otherwise.
This can be
written as
\be
W({\bf s},{\bf D}^2) =
\left\{
\begin{array}{ll}
[1+(-1)^{s_1-s_2}]/2, & s_1 \ge s_2, \\
0,                    & s_1 <   s_2.
\end{array}
\right.
\label{ex3_res}
\ee
As $m=l$ we have ${\bf n}={\bf 0}=(0,0)$, and all
$B_{{\bf 0}}^{(m,m)}
({\bf s}|{\bf D}^m_{\bf j}) = 1$.
The polynomial part is found as
\be
W_{\bf 1}^{\bf 0}({\bf s},{\bf D}^2) =
\frac{\alpha_2}{2(\alpha_1+\alpha_2)}.
\label{ex3_poly}
\ee
In addition to the polynomial part only the
term corresponding to ${\bf j}=(2,2)$ produces nonzero contribution to the
vector partition function, while two other terms with
${\bf j}=(1,2)$ and ${\bf j}=(2,1)$ don't contribute into the result.
The ${\bf j}$-modified matrix ${\bf D}^m_{\bf j}$ reads
($p=2$)
$$
{\bf D}^2_{(2,2)} =
\left(
\begin{array}{ll}
2 \ 2  \\
0 \ 2
\end{array}
\right),
$$
and the corresponding term is found
\be
W_{(2,2)}^{\bf 0}({\bf s},{\bf D}^2_{(2,2)}) =
\frac{\alpha_2}{4(\alpha_1+\alpha_2)}
\sum_{r=0}^1 \Psi_2(s_1-r)\Psi_2(s_2-r) =
\frac{\alpha_2}{2(\alpha_1+\alpha_2)}(-1)^{s_1-s_2}.
\label{ex3_wave}
\ee
Combining (\ref{ex3_poly}) and (\ref{ex3_wave}) we arrive at the expression
$$
%\be
W({\bf s},{\bf D}^2) =
\frac{\alpha_2}{2(\alpha_1+\alpha_2)}[1+(-1)^{s_1-s_2}].
$$
%\label{ex3_final}
%\ee
This result also can be obtained from %formula
(\ref{example0}) applying to it the
recursive relation (\ref{Wvectrecurs}) with ${\bf c}_3=(0,1)$.
It is easy to check that 
the choice $\alpha_1/\alpha_2=0$ corresponds to the first
line in (\ref{ex3_res}), while $\alpha_2/\alpha_1=0$ produces the second one.

%%%%%%%%%%%%%%%%%%%%%%%%%%%%%%%%%%%%%%%%%%%%%%%%%%%%%%%%%%%%%%%%%%%%%%%%%%%%%%%%%
\subsection{Example 4}
Finally consider the case $m=4, l=3$ with
$$
{\bf D}^4 =
\left(
\begin{array}{ll}
2 \ 1 \ 0 \ 0  \\
0 \ 1 \ 1 \ 2  \\
0 \ 0 \ 1 \ 0
\end{array}
\right), \ \ \ \
{\bf s} = (s_1,s_2,s_3).
$$
The straightforward computation using the partial fraction expansion gives the
following result
\be
W({\bf s},{\bf D}^4) =
\left\{
\begin{array}{ll}
\left(1+(-1)^{s_1+s_2+s_3}\right)
\left(\frac{3+(-1)^{s_1}}{8}+\frac{s_1}{4}\right), &
s_2 \ge s_3, \;s_1-s_2+s_3-2 < 0, \\
\left(1+(-1)^{s_1+s_2+s_3}\right)
\left(\frac{3+(-1)^{s_1}}{8}+\frac{s_2-s_3}{4}\right), &
s_2 \ge s_3, \;s_1-s_2+s_3-2 \ge 0, \\
0,                    & s_2 <   s_3.
\end{array}
\right.
\label{ex4_res}
\ee
There are three partial waves
corresponding to
${\bf n}_1=(1,0,0),\;{\bf n}_2=(0,1,0),\;{\bf n}_1=(0,0,1)$.
We drop the computation details and
present the final result in the form
$$
%\be
W({\bf s},{\bf D}^4) =
\frac{\alpha_3[1+(-1)^{s_1+s_2+s_3}]}{8(\alpha_1+\alpha_2)(\alpha_2+\alpha_3)}
%$$
%$$
%\times
[3\alpha_1+4\alpha_2+(-1)^{s_1}(\alpha_1+\alpha_2)
+(-1)^{s_3}(\alpha_2+\alpha_3)+
2\mbf{\alpha}\cdot{\bf s}].
$$
%\label{ex3_final}
%\ee
The comparison of the real part of
the above expression with (\ref{ex4_res}) leads to the following
values of the vector $\mbf{\alpha}$ -- the first line in
(\ref{ex4_res}) corresponds to
$\alpha_2/\alpha_1=0,
\alpha_3/\alpha_1=\pm i$, the second one is obtained
with 
$\alpha_2/\alpha_1=\pm i,
\alpha_3/\alpha_1=1\mp i$, and
the last line
is given by $\alpha_3/\alpha_1=
\alpha_3/\alpha_2=0$.

\subsection{Discussion}
The formula (\ref{Wvectfinal}) provides an exact solution of the
vector restricted partition problem through the vector Bernoulli polynomials of
higher order.
% and it contains free vector parameter $\mbf{\alpha}$.
It also can be written in the form
of a mixture of partial waves (\ref{Wvectfinalpartial})
and the superposition coefficients depend on the vector $\mbf{\alpha}$.
It appears
that proper choice of this vector parameter produces
the solution in one of the system chambers.
The procedure for selection of the vector $\mbf{\alpha}$ will be published
elsewhere.

Knowledge of the vector $\mbf{\alpha}$  for
each chamber gives way to
determine chamber boundaries. The partition function
(as well as its Sylvester waves) in two adjacent
chambers should coincide at their boundary. Thus, setting the
corresponding values to $\mbf{\alpha}$ and equating the
real part of resulting expressions
for $W({\bf s},{\bf D})$ one obtains the linear relations between components of
the vector ${\bf s}$ which determines the location of the boundary.

Thus, the formulas (\ref{Wvectfinal}) and (\ref{Wvectfinalpartial})
appear to be the most general as they
not only give the expression for the vector partition function in every chamber
of the system, but also determine shape of each individual chamber.

The derivation of (\ref{Wvectfinal}) was made in assumption that the
matrix {\bf D} is non-degenerate.
The degeneracy of the matrix related to linear dependence of its
columns doesn't affect the solution of the problem.
It happens only if the matrix rows are linearly dependent.
This dependence leads to effective reduction of the
row number $l$ and
to additional linear conditions imposed on the components of the
vector ${\bf s}$.
These conditions should be satisfied to have nonvanishing value of the
vector partition function.
This way the original problem is reduced to the
non-degenerate case.

\end{document}